\documentclass[12pt]{amsart}

\usepackage{amsthm,amssymb,amsmath,epsfig,graphics,appendix,bm, mathrsfs,bbm}

\usepackage{xcolor}
\usepackage[T1]{fontenc}
\usepackage{stmaryrd} 
\usepackage[left=1in, right=1in, top=1.1in,bottom=1.1in]{geometry}
\usepackage{enumitem}
\usepackage{hyperref}
\hypersetup{
     colorlinks   = true,
     citecolor    = blue,
     linkcolor=blue
}

\usepackage{comment}

\usepackage{csquotes}
\usepackage{graphicx}
\usepackage{pstricks}
\usepackage{lmodern}

\numberwithin{equation}{section}

\newtheorem{thm}{Theorem}[section]
\newtheorem{lemma}[thm]{Lemma}

\newtheorem{corollary}[thm]{Corollary}

\def\E{\mathbb E}
\def\R{\mathbb R}

\DeclareTextFontCommand{\emph}{\em}

\begin{document}

\title[Small deviation of
intersection local time]{Small deviation for
  the mutual intersection local time of Brownian motions}
\author[X. Chen]{Xia Chen}
\address{Department of Mathematics, University of Tennessee,   Knoxville}
\email{xchen@math.utk.edu}

\author[J. Song]{Jian Song}
\address{Research Center for Mathematics and Interdisciplinary Sciences, Shandong University}
\email{txjsong@sdu.edu.cn}

\begin{abstract}
In this note, we establish the bounds

\[ c\varepsilon^{\frac23}\le P\bigg\{\int_0^1\!\!\int_0^1\delta_0(B_s-\tilde{B}_r)dsdr\le \varepsilon \bigg\} \le C \varepsilon^{\frac23},\]
for the mutual intersection local time of two independent 1-dimensional Brownian motions
 $B$ and $\tilde B$.

\end{abstract}

\date{}
\maketitle

\section{Introduction}

Intersection local times (including self-intersection local time and mutual intersection local time) of stochastic processes are a fundamental concept in stochastic 
analysis, which also play an important role in the study of the quantum field theory,  random polymer models, stochastic partial differential
equations and the models of polymers. Our goal in this note is to
establish the small deviation for the mutual intersection local time
$$
\int_0^t\!\int_0^t\delta_0(B_s-\tilde{B}_r)dsdr,\hskip.2in t>0
$$
driven by two independent 1-dimensional Brownian motions $B_s$ and
$\tilde{B}_s$ with $B_0=\tilde{B}_0=0$, where $\delta_0(\cdot)$ is the Dirac delta
function on $\R$ (see, e.g.,   \cite[Chapter 2]{Chen10} for the detailed construction of  intersection local times).
The mutual  intersection local time measures the intensity of path intersection between two
independent Brownian trajectories.
The following is the main result of this paper.

\begin{thm}\label{thm:main-result}
There exists constants $0<c\le C<\infty$ such that for small $\varepsilon>0$,
\begin{equation}\label{e:main-result}
c\varepsilon^{\frac23}\le P\bigg\{\int_0^1\!\!\int_0^1\delta_0(B_s-\tilde{B}_r)dsdr\le \varepsilon \bigg\} \le C \varepsilon^{\frac23}.
\end{equation}
\end{thm}

The small deviation (or, small ball probability) is an interesting subject
for stochastic models. An interested reader is referred to \cite{ls}
for a survey of general development in this direction. 
The goal is to find the  decay rate of the probability
$P\{\|G\|\le\varepsilon\}$ as $\varepsilon\to 0^+$, where $G(\cdot)$ is
a stochastic process or a random field that is properly embedded
in a Banach space endowed with the norm $\|\cdot\|$. In our setting,
$$
\int_0^1\!\!\int_0^1\delta_0(B_s-\tilde{B}_r)dsdr=\int_{-\infty}^\infty L(1,x)\tilde{L}(1,x)dx\buildrel \Delta\over =\|G\|,
$$
where $L(t,x)$ and $\tilde{L}(t,x)$ ($(t,x)\in \R^+\times \R$) are local times of $B$ an $\tilde{B}$, respectively;
and the random field $G(x)=L(1,x)\tilde{L}(1,x)$ ($x\in \R$) is embedded in the space $L^1(\R)$.
Due to technical limitations, the random fields $G(\cdot)$
in the literature of small ball probability are predominantly
Gaussian. Theorem \ref{thm:main-result} represents a rare setting
for the non-Gaussian situation. Another similar example of the non-Gaussian is the small deviation for
the self-intersection local time (see (\ref{self}) below).

In addition, this work is also motivated by a practical matter arising from the area of  stochastic
partial differential equations:
 The integrability
\begin{align}\label{-p}
\E\bigg[\int_0^1\!\!\int_0^1\delta_0(B_s-\tilde{B}_r)dsdr\bigg]^{-p}<\infty \text{ for some } p>0,
\end{align}
plays an important rule when studying the regularity of the probability distribution for the solution of 1-dimensional parabolic Anderson equation (driven by a
space-white Gaussian noise) in \cite{hns11}. For possible quantification of the smoothness of the solution density in the future, it might be desirable to find the critical value of $p$.  The following result is a direct consequence of \eqref{e:main-result}:

\begin{corollary} The negative-moment-integrability (\ref{-p}) holds
if and only if $p<\frac23$.
\end{corollary}

The large and small deviations are two research subjects concerning the probabilities of the rare events
representing two extreme behaviors
of the stochastic system: the random variable takes
 unusually large values and the random variable takes unusually small
values. For the purpose of comparison in terms of large and small deviations, we introduce the self-intersection
local time
$$
\int_0^t\!\int_0^t\delta_0(B_s-B_r)dsdr,\hskip.2in t>0.
$$
By Theorem 1.1 in \cite{cl04} ($m=1$, $p=2$ for self-intersection, and $m=2$,
$p=1$ for mutual intersection), we have the large deviations for both self and mutual
intersections:
$$
\lim_{\varepsilon\to 0^+}{\varepsilon^2}\log P\bigg\{\int_0^1\!\int_0^1\delta_0(B_s-B_r)dsdr
\ge \varepsilon^{-1}\bigg\}=-{3\over 2};
$$
$$
\lim_{\varepsilon\to 0^+}{\varepsilon^2}
\log P\bigg\{\int_0^1\!\int_0^1\delta_0(B_s-\tilde{B}_r)dsdr
\ge \varepsilon^{-1}\bigg\}=-3.
$$
It is noticeable that the large deviations for two different types
of intersections take
the forms close to each other. The story behind is the following deterministic relation:
$$
\begin{aligned}
 &\int_0^1\!\int_0^1\delta_0(B_s-\tilde{B}_r)dsdr
=\int_{-\infty}^\infty L(1,x)\tilde{L}(1,x)dx
\le{1\over 2}\int_{-\infty}^\infty L^2(1,x)dx+{1\over 2}\int_{-\infty}^\infty
   \tilde{L}^2(1,x)dx\\
  &={1\over 2}\int_0^1\!\int_0^1\delta_0(B_s-B_r)dsdr
    +{1\over 2}\int_0^1\!\int_0^1\delta_0(\tilde{B}_s-\tilde{B}_r)dsdr
    \end{aligned}
$$
where $L(t,x)$ and $\tilde{L}(t,x)$ are, respectively,  the local times of
$B_t$ and $\tilde{B}_t$.
In the game of large deviations, the
Brownian trajectories that make mutual intersection large are those
that make $L(1,\cdot)$ and $\tilde{L}(1,\cdot)$ close to each
other so that ``$\le$''
is replaced by ``$\approx$'' in the above deterministic relation.

The small ball probability for the self-intersection local time
(\cite{hhk03}) takes the form
\begin{align}\label{self}
\lim_{\varepsilon\to 0^+} \varepsilon^2\log P\left\{\int_0^1\!\!\int_0^1\delta_0(B_s- B_r)dsdr\le \varepsilon \right\}=-c
\end{align}
with constant  $c>0$.  A striking difference
between (\ref{self}) and (\ref{thm:main-result})
lies in the disparity between exponential decay and power decay.
The typical paths that the Brownian motion takes for maximizing
self-avoiding  are the ones that are close to straight lines, i.e., $B_t\approx c\epsilon^{-1}t, t\in[0,1].$ For the mutual intersection local time, 
the small deviation means a totally different game: Two independent Brownian paths
tend to avoid meeting each other right after they  set out from the same starting point.

In the following we outline
the idea in our proof: By the Brownian scaling property and the relation
$\delta_0(\lambda x)=\lambda^{-1}\delta_0(x)$ for any $\lambda >0$, we have
\begin{align}\label{scalling}
\int_0^t\!\int_0^t\delta_0(B_s-\tilde{B}_r)dsdr\buildrel d\over =t^{3/2}\int_0^1\!\int_0^1\delta_0(B_s-\tilde{B}_r)dsdr, \hskip.2in \forall t>0,
\end{align}
where ``$\buildrel d\over =$'' means equality in distribution. Therefore, the proof of (\ref{thm:main-result}) is reduced to establishing the following bound
\begin{align}\label{bounds}
c{1\over t}\le P\bigg\{\int_0^t\!\!\int_0^t\delta_0(B_s-\tilde{B}_r)dsdr\le a \bigg\} \le C {1\over t}
\end{align}
for large $t$ and for some fixed constants $c, C, a>0$ which are independent of $t$.

In the proof of the lower bound,  we separate the Brownian paths $B[1,t]$ and $\tilde{B}[1,t]$ 
into two disjoint and distanced
half-lines in a way that leads to $B[0,t]\cap \tilde{B}[1,t]=\emptyset$ and $B[1,t]\cap \tilde{B}[0,t]=\emptyset$.
Consequently, the only unavoidable intersections are the ones happening 
on the time square $[0,1]^2$.

In other words, the strategy lies in the comparison
between the probability of ``low  intersection'' and the probability of ``no intersection''. 
Indeed, the steps taken in  (\ref{e:2.2}) below morally suggests the relation
\begin{align}\label{inter}
P\{\hbox{low intersection on $[0,t]^2$}\}\ge cP\{\hbox{no intersection on $[1,t]^2$}\}.
\end{align}
Since the probability on the right-hand side is not hard to compute
(only for  1-dimensional Brownian motions, of course),
 the proof of the lower bound (see Section \ref{sec:lower}) is the easier one, in comparison with the proof of the upper bound (see Section \ref{sec:upper}).
As for the upper bound, the real challenge is to fill the gap between ``low intersection'' and ``no intersection'' suggested by~(\ref{inter}).

We end the discussion with a remark on the possible multi-dimensional
extensions in the future. It is well-known
that the mutual intersection local time between two indpendent
$d$-dimensional Brownian motions exists for $d=1,2,3$.
Therefore, it  makes sense to raise the same question in these dimensions. In view of the heuristic comparison (\ref{inter}),\footnote{
We believe it remains morally true in multi-dimensional settings.}
it brings the famous theorem of intersection exponent by Lawler, Schramm and Werner (Theorem 1.1,~\cite{lsw}) in which it claims that when $d=2$, the right-hand
side of (\ref{inter}) is of the order  $t^{-(5/8)+o(1)}$ as $t\to\infty$. When $d=2$, noting that
$$
\int_0^t\!\int_0^t\delta_0(B_s-\tilde{B}_r)dsdr\buildrel d\over =t\int_0^1\!\int_0^1\delta_0(B_s-\tilde{B}_r)dsdr,\hskip.2in \forall t>0,
$$
 we conjecture that
\begin{align}\label{d=2}
P\bigg\{\int_0^1\!\int_0^1\delta_0(B_s-\tilde{B}_r)dsdr\le\varepsilon\bigg\}=\varepsilon^{(5/8)+o(1)},\hskip.2in (\varepsilon\to 0^+).
\end{align}

\medskip

\noindent{\bf Notations.} In the proofs, the Brownian motions $B$ and $\tilde B$ may start from places other than the origin.  We use $P_{(x,\tilde{x})}(\cdot)$ for the probability with $B_0=x$ and $\tilde{B}_0=\tilde{x}$, and we take the convention $P(\cdot) \buildrel \Delta\over = P_{(0,0)}(\cdot)$. We also denote $P_x(A) \buildrel \Delta\over = P_{(x,\tilde x)} (A)$ and $P(A)  \buildrel \Delta\over = P_0(A)$ for $A\in\sigma(B_s, 0\le s<\infty)$ and similarly $P_{\tilde x}(A) \buildrel \Delta\over = P_{(x,\tilde x)} (A)$ and $P(A) \buildrel \Delta\over = P_0(A)$ for $A\in\sigma(\tilde B_s, 0\le s<\infty)$. The letters $c$ and $C$ denote generic positive constants that may vary in different places.

\section{Lower bound}\label{sec:lower}

In this section, all we need is to prove the lower bound in (\ref{bounds}) with $a=1$, i.e.,
\begin{equation}\label{e:lb1}
P\bigg\{\int_0^t\!\!\int_0^t\delta_0(B_s-\tilde{B}_r)dsdr\le 1\bigg\}\ge \frac ct
\end{equation}
for large $t$. First, we claim that for $t>1$,
\begin{equation}\label{e:2.2}
\begin{aligned}
&P\bigg\{\int_0^t\!\!\int_0^t\delta_0(B_s-\tilde{B}_r)dsdr\le 1\bigg\}\\
&\ge
  P\bigg\{\int_0^1\!\!\int_0^1\delta_0(B_s-\tilde{B}_r)dsdr\le 1,\hskip.05in \min_{0\le s\le 1}B_s\ge -1,
  \hskip.05in \max_{0\le s\le 1}\tilde{B}_s\le 1,\\
&\quad \quad 
 \min_{1\le s\le t}B_s> 1,
\hskip.05in \max_{1\le s\le t}\tilde{B}_s<-1\bigg\}\\
&=\E \Bigg[\mathbf{\mathbf 1}_A  \cdot P_{B_1}\bigg\{ \min_{0\le s\le t-1}B_s
  >1\bigg\}
  P_{\tilde{B}_1}\bigg\{
  \max_{0\le s\le t-1}\tilde{B}_s< -1\bigg\}\Bigg],
  \end{aligned}
\end{equation}
where 
\[A=\bigg\{\int_0^1\!\!\int_0^1\delta_0(B_s-\tilde{B}_r)dsdr\le 1,\hskip.05in \min_{0\le s\le 1}B_s\ge -1,
  \hskip.05in \max_{0\le s\le 1}\tilde{B}_s\le 1\bigg\},\]
and the last step follows from the Markov property and the independence between $B$ and~$\tilde{B}$.

Indeed, the first step is justified by the fact that
 on the event  
 $$\Big\{\min_{0\le s\le 1}B_s\ge -1,
  \hskip.05in \max_{0\le s\le 1}\tilde{B}_s\le 1,\hskip.05in
 \min_{1\le s\le t}B_s> 1,
\hskip.05in \max_{1\le s\le t}\tilde{B}_s<-1\Big\},
$$
$B[0,t]\cap \tilde{B}[1,t]=\emptyset$
and $B[1,t]\cap \tilde{B}[0,t]=\emptyset$ which leads to
$$
\int_0^t\!\!\int_0^t\delta_0(B_s-\tilde{B}_r)dsdr
=\int_0^1\!\!\int_0^1\delta_0(B_s-\tilde{B}_r)dsdr.
$$

On $\{B_1\ge 2\}$, by (\ref{reflection-1}) in  Lemma \ref{lem:reflection} we have
\begin{equation}\label{e:2.3}
\begin{aligned}
&P_{B_1}\bigg\{ \min_{0\le s\le t-1}B_s>1\bigg\}
\ge P\bigg\{1+\min_{0\le s\le t-1}B_s
  >0\bigg\}\\
&  =P\bigg\{\vert B_{t-1}\vert\le 1\bigg\}\sim {1\over \sqrt{2\pi t}}\hskip0.2in (t\to\infty),
  \end{aligned}
\end{equation}
and similarly, using (\ref{reflection-2}) in  Lemma \ref{lem:reflection}  we conclude that
on $\{\tilde{B}_1\le -2\}$,
\begin{equation}\label{e:2.4}
P_{\tilde{B}_1}\bigg\{
  \max_{0\le s\le t-1}\tilde{B}_s< -1\bigg\}\ge P\{\vert B_{t-1}\vert\le 1\}
  \sim {1\over \sqrt{2\pi t}} \hskip0.2in (t\to\infty). 
\end{equation}

Combining \eqref{e:2.2}, \eqref{e:2.3} and \eqref{e:2.4}, we get that when $t$ is sufficiently large, there exists a constant $c>0$ such that
$$
\begin{aligned}
&P\bigg\{\int_0^t\!\!\int_0^t\delta_0(B_s-\tilde{B}_r)dsdr\le 1\bigg\}\\
&\ge \frac ct \, P\bigg\{\int_0^1\!\!\int_0^1\delta_0(B_s-\tilde{B}_r)dsdr\le 1,\hskip.05in 
  \min_{0\le s\le 1}B_s\ge -1,\hskip.05in  B_1\ge 2,\hskip.05in\max_{0\le s\le 1}\tilde{B}_s\le 1,\hskip.05in \tilde{B}_{1}\le -2
\bigg\},
\end{aligned}
$$
where the probability on the right-hand side is clearly positive. The desired lower bound~\eqref{e:lb1} is proved.

\begin{lemma}\label{lem:reflection}
For any $x>0$ and $t>0$, we have
\begin{align}\label{reflection-1}
P\{x+B_s>0;\hskip.05in \forall s\le t\}
=P\{\vert B_t\vert\le x\},
\end{align}
and 
\begin{align}\label{reflection-2}
P\{-x+B_s<0;\hskip.05in \forall s\le t\}
=P\{\vert B_t\vert\le  x\}.
\end{align}
\end{lemma}

\proof Denote $\chi_x\buildrel \Delta\over=\inf\{s\ge 0;\hskip.05in B_s=x\}$.  Then, 
$$
P\{x+B_s>0;\hskip.05in \forall s\le t\}=P\{\chi_{-x}\ge t\}=P\{\chi_x\ge t\}
=P\Big\{\max_{0\le s\le t}B_s\le x\Big\}=P\{\vert B_t\vert\le x\},
$$
where the last equality follows from the reflection principle $\max\limits_{0\le s\le t}B_s\buildrel d\over =\vert B_t\vert$. This yields (\ref{reflection-1}). Equation (\ref{reflection-2}) can be obtained by replacing $B$ by $-B$ in (\ref{reflection-1}).

\section{Upper bound}\label{sec:upper}

In this section, we prove the upper bound of (\ref{bounds}) which can be slightly rephrased as
\begin{equation}\label{e:ub}
P\bigg\{\int_0^t\!\!\int_0^t\delta_0(B_s-\tilde{B}_r)dsdr\le a\bigg\}=O\Big(\frac 1t\Big)\hskip.2in (t\to\infty)
\end{equation}
for some constant $a>0$. We shall specify the constant $a>0$ later.

As we mentioned earlier, the central piece of the proof is to compare the ``low intersection'' probability in (\ref{e:ub}) with
the probability of ``no intersection'' in a direction opposite to the one in (\ref{inter}).
More precisely, we compare   the probability of low intersection for $B$ and $\tilde{B}$ with
the probability of no intersection for 
two independent simple random walks. We refer to
(\ref{Lap-1}) and (\ref{e:2k1})  below for the implementation of
 this idea.

To this end,  we generate two independent simple random walks
$\{S_k, k\ge 0\}$ and $\{\tilde S_k, k\ge 0\}$ by Brownian motions 
in a way that the trajectory $\{S_k;\hskip.05in 0\le k\le n\}$ is intimately close to the Brownian path
$\{B_s;\hskip.05in 0\le t\le\tau_n\}$
up to  the stopping time $\tau_n$ that is given below (same thing is expected between $\tilde{S}_k$ and
$\tilde{B}_t$):

We first define
a sequence of stopping times $\{\tau_k\}_{k\ge 0}$ related to the Brownian motion $B_t$:
\begin{equation}\label{e:tau}
\tau_0=0; ~ \tau_k=\inf\{t\ge \tau_{k-1};\hskip.05in \vert B_t-B_{\tau_{k-1}}\vert =1\}  \text{ for } k=1, 2, \dots
\end{equation}
Then $\{S_k=B_{\tau_k}, k\ge 0\}$) is a simple random walk with $S_0=B_0$. The stopping times 
$\{\tilde{\tau}_k, k\ge 0\}$ and the simple random walk $\{\tilde{S}_k, k\ge 0\}$ are generated from $\{\tilde B_t,t\ge 0\}$ in the same way.

Clearly, $\{\tau_k-\tau_{k-1}\}_{k\ge 1}$ is a sequence of i.i.d.  random variables with the same distribution as $\tau_1$. Here we point out
that the distribution  of $\tau_1$ does not depend on 
the starting point of $B$.\footnote{We mention this as we allow $B_0\not=0$
  and $\tilde{B}_0\not =0$ in later discussion.}
The exact distribution of $\tau_1$
is known (p.342, \cite{Feller}). An instructive link  (\cite{ls}) to the
integrability of $\tau_1$ is given as follows:
By the definition of $\tau_1$,
$$
P\{\tau_1\ge t\}=P\Big\{\max_{s\le t}\vert B_s\vert\le 1\Big\}
=P\Big\{\max_{s\le 1}\vert B_s\vert\le {1\over\sqrt{t}}\Big\}
=\exp\Big\{-\Big({\pi^2\over 8}+o(1)\Big)t\Big\} \hskip0.15in (t\to\infty),
$$
where the last step is the classic result on the small
ball probability for Brownian motions (see, e.g., (1.3), \cite{ls}).
In particular, $\E e^{\theta\tau_1}<\infty$
for $\theta<\pi^2/8$. By a standard application of Chebyshev's inequality to the sum of independent random variables, we get
\begin{align}\label{Chebyshev}
  P\big\{\vert \tau_n-n\E\tau_1\vert\ge n\delta\big\}
  \le {1\over n\delta^2} \text{Var}(\tau_1).
\end{align}

The mutual intersection local time $Q_n$ of the random walks $S$ and $\tilde S$ is given by
\begin{equation}\label{e:Qn}
Q_n\buildrel \Delta\over= \sum_{j,k=1}^n{\mathbf 1}_{\{S_j=\tilde{S}_k\}}.
\end{equation}

For $n\in\mathbb N$,  
\begin{equation}\label{minor}
\begin{aligned}
  &\int_0^{\tau_{n+1}}\!\!\int_0^{\tilde{\tau}_{n+1}}
    \delta_0(B_s-\tilde{B}_r)dsdr=\sum_{j,k=0}^{n}\int_{\tau_{j}}^{\tau_{j+1}}\!\!
    \int_{\tilde{\tau}_{k}}^{\tilde{\tau}_{k+1}}\delta_0(B_s-\tilde{B}_r)dsdr\\
&     \ge  \sum_{j,k=1}^{n}\int_{\tau_{j}}^{\tau_{j+1}}\!\!
    \int_{\tilde{\tau}_{k}}^{\tilde{\tau}_{k+1}}\delta_0(B_s-\tilde{B}_r)dsdr\ge \sum_{j,k=1}^n{\mathbf 1}_{\{S_j=\tilde{S}_k\}}\int_{\tau_{j}}^{\tau_{j+1}}\!\!
    \int_{\tilde{\tau}_{k}}^{\tilde{\tau}_{k+1}}\delta_0(B_s-\tilde{B}_r)dsdr\\
  &=\sum_{j,k=1}^n{\mathbf 1}_{\{S_j=\tilde{S}_k\}}\xi_{j,k}=H_n\hskip.2in \hbox{(say)},
    \end{aligned}
    \end{equation}
where we denote
\[\xi_{j,k}\buildrel \Delta\over= \int_{\tau_{j}}^{\tau_{j+1}}\!\!
    \int_{\tilde{\tau}_{k}}^{\tilde{\tau}_{k+1}}\delta_0(B_s-\tilde{B}_r)dsdr \quad \text{ for } j, k =1, 2, \dots\]

For a technical reason, we allow the Brownian motions to start form somewhere different from $0$ and use
$P_{(x,\tilde{x})}(\cdot)$ for the probability with $B_0=x$ and $\tilde{B}_0=\tilde{x}$.
We follow the convention $P(\cdot)=P_{(0,0)}(\cdot)$. Recall the standard notation used in
the textbook on Markov process:
$$
P_\mu(\cdot)=\int_{\R^2}\mu(dx,d\tilde{x})P_{(x,\tilde{x})}(\cdot)
$$
for any probability measure $\mu$ on $\R^2$. In the following discussion, we set
$$
\mu={1\over 4}\Big\{\delta_{(-1, -1)}+\delta_{(1,1)}+\delta_{(-1,1)}+\delta_{(1,-1)}\Big\}.
$$
Despite the fact that
$(0,0)$ is not listed as a starting point of the
2-dimensional Brownian motion $(B_t, \tilde{B}_t)$ under the distribution $\mu$, noting that
\begin{align}\label{notation-1}
  \left\{\begin{array}{ll}\displaystyle
    P_{(-1,-1)}\{H_n\le a\}=P_{(1,1)}\{H_n\le a\}=P\{H_n\le a\},\\\\
    \displaystyle P_{(-1,1)}\{H_n\le a\}=P_{(1,-1)}\{H_n\le a\}\hskip.2in
  n=1,2,\cdots,\end{array}\right. 
\end{align}
we have
\begin{align}\label{notation-2}
P_\mu\{H_n\le a\}={1\over 2}P\{H_n\le a\}+{1\over 2}P_{(1,-1)}\{H_n\le a\}.
\end{align}
Here we recall the notation from \eqref{minor},
$$
H_n=\sum_{j,k=1}^n{\mathbf 1}_{\{S_j=\tilde{S}_k\}}\xi_{j,k}.
$$
The proof of (\ref{e:ub}) can be reduced to  establishing the bound
  \begin{align}\label{e:Hn}
  P_\mu\{H_n\le a\}=O\Big({1\over n}\Big)\hskip.2in (n\to\infty).
  \end{align}

  Indeed, by (\ref{notation-2}) we have $P\{H_n\le a\}
  \le 2P_\mu\{H_n\le a\}$. Therefore, (\ref{minor}) leads to
    $$
    P\bigg\{\int_0^{\tau_{n+1}}\!\!\int_0^{\tilde{\tau}_{n+1}}\delta_0(B_s-\tilde{B}_r)dsdr\le a\bigg\}=O\Big({1\over n}\Big)\hskip.2in (n\to\infty).
    $$
    Consequently, for any $0<\delta<1$,
    $$
    \begin{aligned}
    &P\bigg\{\int_0^t\!\!\int_0^t \delta_0(B_s-\tilde{B}_r)dsdr\le a\bigg\}\\
    &\le
    P\bigg\{\int_0^{\tau_{[(1-\delta)t/\E\tau_1]+1}}\!\!\int_0^{\tilde{\tau}_{ [(1-\delta)t/\E\tau_1]+1}}\delta_0(B_s-\tilde{B}_r)dsdr\le a\bigg\}
    +2P\Big\{\tau_{ [(1-\delta)t/\E\tau_1]+1}\ge t\Big\}\\
    &\le O\Big({1\over t}\Big)+{C\over t}=O\Big({1\over t}\Big) \hskip.2in (t\to\infty),
    \end{aligned}
    $$
  where the second inequality partially follows from (\ref{Chebyshev}).

  \medskip
  
The remaining of the section is devoted to the proof of (\ref{e:Hn}). Set the random set on $\mathbb Z_+^2$:
$$
\Lambda(A)\buildrel \Delta\over=\{(j,k)\in A;\hskip.05in S_j=\tilde{S}_k\},
\hskip.2in A\subset \mathbb Z_+^2.
$$
Notice that  $\#(\Lambda([1, n]^2)=Q_n$ where $Q_n$ is the number of intersections given in \eqref{e:Qn}. We introduce the following stopping time 
\begin{equation}\label{e:sigma}
\sigma\buildrel \Delta\over=\min\left\{n\ge 1;\hskip.05in \Lambda([1,n]^2)\neq \emptyset\right\}=\min\big\{n\ge 1; Q_n>0\big\}.
\end{equation}

 By the law of total probability, we have
\begin{equation}\label{e:bd-Hn}
\begin{aligned}
P_\mu\{H_n\le a\}&=\sum_{l=1}^nP_\mu\{\sigma=l,\hskip.05in H_n\le a\}+P_\mu\{\sigma>n\}\\
&=\sum_{l=1}^nP_\mu\{\sigma=l,\hskip.05in H_n\le a\}+P_\mu\{Q_n=0\}.
\end{aligned}
\end{equation}
On the event $\{\sigma=l\}$, the intersection $\{S_j=\tilde{S}_k\}$ happens either on $[1,l]\times\{l\}$ or on $\{l\}\times[1,l]$,
or on both. An  intersection $\{S_j=\tilde{S}_k\}$ is called an early intersection, if there is no $(j_1, k_1)\not =(j,k)$
such that $j_1\le j$, $k_1\le k$, and $S_{j_1}=\tilde{S}_{k_1}$. In general, there may be multiple early intersections
due to the  multi-dimension of time.
On $\{\sigma=l\}$, early intersections happen on $[1,l]\times\{l\}$ or on $\{l\}\times[1,l]$, or on both, i.e.,
$$
\begin{aligned}
\{\sigma=l\}=&\left\{\Lambda([0,l-1]^2)=\emptyset, \hskip.05in \hbox{early intersection happens on $[1,l]\times\{l\}$}\right\}\\
&\bigcup
\left\{\Lambda([0,l-1]^2)=\emptyset, \hskip.05in \hbox{early intersection happens on $\{l\}\times [1, l]$}\right\}.
\end{aligned}
$$
Therefore, for $l=1, \dots, n$,
\begin{equation}\label{e:l-Hn}
\begin{aligned}
&P_\mu\{\sigma=l,\hskip.05in H_n\le a\}\\
&\le P_\mu\left\{\Lambda([0,l-1]^2)=\emptyset, \hskip.05in \hbox{early intersection happens on $[1,l]\times\{l\}, \hskip.05in H_n\le a$}\right\}\\
&\quad +P_\mu\left\{\Lambda([0,l-1]^2)=\emptyset, \hskip.05in \hbox{early intersection happens on $\{l\}\times [1,l], \hskip.05in H_n\le a$}\right\}\\
&=2P_\mu\left\{\Lambda([0,l-1]^2)=\emptyset, \hskip.05in \hbox{early intersection happens on $[1,l]\times\{l\}, \hskip.05in H_n\le a$}\right\}.
\end{aligned}
\end{equation}

On the event $\big\{\Lambda([0,l-1]^2)=\emptyset, \hskip.05in \hbox{early intersection happens on $[1,l]\times\{l\}$}\big\}$, there is
a unique early intersection on $[1,l]\times\{l\}$
happening at $(\rho, l)$ for some $1\le \rho\le l$, described by the event
\begin{equation}\label{e:F-rho-l}
F_{\rho,l}\buildrel \Delta\over=\Big\{\Lambda\Big(\big([1,\rho]\times[1,l]\big)\setminus (\rho, l)\Big)=\emptyset,
\hskip.05in S_\rho=\tilde{S}_l\Big\}, \hskip.2in \rho=1,\cdots, l.
\end{equation}
Thus, we have
$$
\left\{\Lambda([0,l-1]^2)=\emptyset, \hskip.05in \hbox{early intersection happens on $[1,l]\times\{l\}$}\right\}
\subset \bigcup_{\rho=1}^lF_{\rho,l}\,.
$$
Here we mention that, when $\rho =l$, the fact that ``$S_l=\tilde{S}_l$'' is an early intersection requires, in addition to 
``$S_1\not=\tilde{S}_l,\cdots, S_{l-1}\not =\tilde{S}_l$'',
that
``$S_l\not =\tilde{S}_1,\cdots, S_l\not =\tilde{S}_{l-1}$'' (and therefore $\Lambda\left(\left([1,l]\times[1,l]\right)\setminus (l, l)\right)=\emptyset$). Thus,
$$
\begin{aligned}
&P_\mu\left\{\Lambda([1,l-1]^2)=\emptyset, \hskip.05in \hbox{early intersection happens on $[1,l]\times\{l\}$}, \hskip0.05in H_n\le a\right\}\\
&\le \sum_{\rho=1}^lP_\mu\big\{F_{\rho,l}\cap\{H_n\le a\}\big\}=\sum_{\rho=1}^lP_\mu\big\{F_{\rho,l}\cap\{\xi_{\rho, l}\le a\}\cap\{H_n\le a\}\big\}.
\end{aligned}
$$
Combining this with \eqref{e:l-Hn} yields
\begin{equation}\label{e:l-Hn1}
P_\mu\{\sigma=l,\hskip.05in H_n\le a\}\le 2\sum_{\rho=1}^lP_\mu\big\{F_{\rho,l}\cap\{\xi_{\rho, l}\le a\}\cap\{H_n\le a\}\big\}, \hskip0.2in l=1, \dots, n.
\end{equation}

When $1\le l\le n-2$,   we have  on $F_{\rho,l}$, 
\begin{equation}\label{e:H-H}
\begin{aligned}
&H_n\ge\sum_{j=\rho+2}^n\sum_{k=l+2}^n{\mathbf 1}_{\{S_j=\tilde{S}_k\}}\xi_{j,k}\\
&=\sum_{j=1}^{n-\rho-1}\sum_{k=1}^{n-l-1}{\mathbf 1}_{\{(S_{\rho+1+j}-S_{\rho+1})+(S_{\rho+1}-S_\rho)=(\tilde{S}_{l+1+k}-\tilde{S}_{l+1})+(\tilde{S}_{l+1}-\tilde S_l)\}}
\xi_{\rho+1+j, l+1+k}\\
&\ge\sum_{j,k=1}^{n-l-1}{\mathbf 1}_{\{(S_{\rho+1+j}-S_{\rho+1})+(S_{\rho+1}-S_\rho)=(\tilde{S}_{l+1+k}-\tilde{S}_{l+1})+(\tilde{S}_{l+1}-\tilde S_l)\}}
\xi_{\rho+1+j, l+1+k}\\
&=\bar H_{n-l-1} \hskip.2in \hbox{(say)}, 
\end{aligned}
\end{equation}
where the equality is due to $S_\rho=\tilde S_l$  on $F_{\rho,l}$.
Note that
$$
\begin{aligned}
&\xi_{\rho+1+j, l+1+k}=\int_{\tau_{\rho+1+j}}^{\tau_{\rho+2+j}}\!\int_{\tilde{\tau}_{l+1+k}}^{\tilde{\tau}_{l+2+k}}
\delta_0(B_s-\tilde{B}_r)dsdr\\
&=\int_{\tau_{\rho+1+j}-\tau_{\rho+1}}^{\tau_{\rho+2+j}-\tau_{\rho+1}}\!\int_{\tilde{\tau}_{l+1+k}-\tilde{\tau}_{l+1}}^{\tilde{\tau}_{l+2+k}-\tilde{\tau}_{l+1}}
\delta_0\Big(B_{\tau_{\rho+1}+s}-\tilde{B}_{\tilde \tau_{l+1}+r}\Big)dsdr\\
&=\int_{\tau_{\rho+1+j}-\tau_{\rho+1}}^{\tau_{\rho+2+j}-\tau_{\rho+1}}\!\int_{\tilde{\tau}_{l+1+k}-\tilde{\tau}_{l+1}}^{\tilde{\tau}_{l+2+k}-\tilde{\tau}_{l+1}}
\delta_0\Big(\big((B_{\tau_{\rho+1}+s}-B_{\tau_{\rho+1}})+(S_{\rho+1}-S_{\rho})\big)\\
&
\quad -\big((\tilde{B}_{\tilde \tau_{l+1}+r}-\tilde{B}_{\tilde \tau_{l+1}})+(\tilde{S}_{l+1}-\tilde{S}_l)\big)\Big)dsdr.
\end{aligned}
$$
Thus, given $(S_{\rho+1}-S_\rho, \tilde S_{l+1}-\tilde S_l)=(z_1, z_2)$, the distribution of $\bar H_{n-l-1}$  is the same as that of $H_{n-l-1}$ under $P_{(z_1, z_2)}$,
and is independent of $\{(B_s, \tilde{B}_r); \hskip.005in 0\le s\le \tau_{\rho+1}\hskip.05in \hbox{and}\hskip.05in 0\le r\le \tilde{\tau}_{l+1}\}$, and is therefore
independent of $F_{\rho,l}$ and $\xi_{\rho,l}$.  By (\ref{e:H-H}) and the Markov property
\begin{equation}\label{e:F-xi-H}
\begin{aligned}
&P_\mu\left\{F_{\rho,l}\cap\{\xi_{\rho, l}\le a\}\cap\{H_n\le a\}\right\}\le P_\mu\left\{F_{\rho,l}\cap\{\xi_{\rho, l}\le a\}\cap\{\bar H_{n-l-1}\le a\}\right\}\\
&=\E_\mu\left[{\mathbf 1}_{F_{\rho, l}}{\mathbf 1}_{\{\xi_{\rho, l}\le a\}}P_{(S_{\rho+1}-S_{\rho}, \tilde{S}_{l+1}-\tilde{S}_l)}\{H_{n-l-1}\le a\}\right].
\end{aligned}
\end{equation}
By the fact that $S_{\rho+1}-S_{\rho}=\pm 1$ and $\tilde{S}_{l+1}-\tilde{S}_l=\pm 1$, and by (\ref{notation-1}),
the right-hand side of \eqref{e:F-xi-H} is bounded from above by
\begin{equation}\label{e:F-xi-H1} \begin{aligned}
&P_{(1,1)}\{H_{n-l-1}\le a\}\E_\mu\left[{\mathbf 1}_{F_{\rho, l}}{\mathbf 1}_{\{\xi_{\rho, l}\le a\}}{\mathbf 1}_{\{S_{\rho+1}-S_{\rho}=\tilde{S}_{l+1}-\tilde{S}_l)}\right]\\
&+P_{(1,-1)}\{H_{n-l-1}\le a\} \E_\mu\left[{\mathbf 1}_{F_{\rho, l}}{\mathbf 1}_{\{\xi_{\rho, l}\le a\}}{\mathbf 1}_{\{S_{\rho+1}-S_{\rho}\not =\tilde{S}_{l+1}-\tilde{S}_l)}\right]\\
&\le \left(P_{(1,1)}\{H_{n-l-1}\le a \}+P_{(1,-1)}\{H_{n-l-1}\le a\}\right)P_\mu\left\{F_{\rho, l}\cap\{\xi_{\rho,l}\le a\}\right\}\\
&=2P_\mu\{H_{n-l-1}\le a\}P_\mu\left\{F_{\rho, l}\cap\{\xi_{\rho,l}\le a\}\right\}.
\end{aligned}
\end{equation}
On $F_{\rho, l}$,
$$
\xi_{\rho, l}=\int_0^{\tau_{\rho+1}-\tau_\rho}\!\int_0^{\tilde{\tau}_{l+1}-\tilde{\tau}_l}\delta_0\Big((B_{\tau_\rho +s}-B_{\tau_\rho})-(\tilde{B}_{\tilde{\tau}_l+r}-\tilde{B}_{\tilde{\tau}_l})\Big)
dsdr
$$
and the right-hand side has the same distribution as $\xi_{1,1}$ (under the law $P(\cdot)$), and is independent of $\big\{(B_s, \tilde{B}_r); \hskip.05in s\le \tau_{\rho}\hskip.05in \hbox{and}\hskip.05in r\le \tilde{\tau}_{l}\big\}$, and
therefore  of $F_{\rho,l}$.
Hence,
\begin{equation}\label{e:F-xi}
  P_\mu\left\{F_{\rho, l}\cap\{\xi_{\rho,l}\le a\}\right\}
  =P_\mu\{F_{\rho, l}\}P\{\xi_{1,1}\le a\}.
\end{equation}

In summary, combining \eqref{e:l-Hn1},  \eqref{e:F-xi-H}, \eqref{e:F-xi-H1} and \eqref{e:F-xi}, we get  for $1\le l\le n-2$,
\begin{equation}\label{e:l-Hn2}
P_\mu\{\sigma=l,\hskip.05in H_n\le a\}\le 4P\{\xi_{1,1}\le a\}\bigg(\sum_{\rho=1}^l P_\mu\{F_{\rho,l}\}\bigg)P_\mu\{H_{n-l-1}\le a\}.
\end{equation}

When $l\in\{n-1,n\}$, by \eqref{e:l-Hn1} and \eqref{e:F-xi} we have the bound
\begin{equation}\label{e:l-Hn3}
P_\mu\{\sigma=l,\hskip.05in H_n\le a\}\le 2P\{\xi_{1,1}\le a\}\bigg(\sum_{\rho=1}^{l} P_\mu\{F_{\rho,l}\}\bigg).
\end{equation}

Now, combining \eqref{e:bd-Hn},  \eqref{e:l-Hn2}, and \eqref{e:l-Hn3}, we get for $n\ge 3$,
\begin{equation}\label{e:key-estimate}
\begin{aligned}
&P_\mu\{H_n\le a\}\\
& \le 4P\{\xi_{1,1}\le a\}\sum_{l=1}^{n-2}\bigg(\sum_{\rho=1}^lP_\mu\{F_{\rho,l}\}\bigg)P_\mu\{H_{n-l-1}\le a\}\\
&\quad +2P\{\xi_{1,1}\le a\}\bigg(\sum_{l\in\{n-1, n\}}\sum_{\rho=1}^{l}P_\mu\{F_{\rho,l}\}\bigg)+P_\mu\{Q_n=0\}.
\end{aligned}
\end{equation}

Therefore, for any $\theta\in(0,1)$,
\begin{equation}\label{e:est1}
\begin{aligned}
&\sum_{n=3}^\infty \theta^nP_\mu\{H_n\le a\}\\
&\le \sum_{n=3}^\infty \theta^nP_\mu\{Q_n=0\}+4P\{\xi_{1,1}\le a\}\sum_{n=3}^\infty\theta^n
\sum_{l=1}^{n-2}\bigg(\sum_{\rho=1}^lP_\mu\{F_{\rho,l}\}\bigg)P_\mu\{H_{n-l-1}\le a\}\\
&\quad +2P\{\xi_{1,1}\le a\}\sum_{n=3}^\infty\theta^n
\bigg(\sum_{\rho=1}^{n-1}P_\mu\{F_{\rho,n-1}\}\bigg)+2P\{\xi_{1,1}\le a\}\sum_{n=3}^\infty\theta^n\bigg(\sum_{\rho=1}^{n}P_\mu\{F_{\rho,n}\}\bigg).
\end{aligned}
\end{equation}
Notice that
\begin{align*}
&\sum_{n=3}^\infty\theta^n\sum_{l=1}^{n-2}\bigg(\sum_{\rho=1}^lP_\mu\{F_{\rho,l}\}\bigg)P_\mu\{H_{n-l-1}\le a\}\\
& = \theta \sum_{n=2}^\infty \sum_{l=1}^{n-1}\bigg(\theta^l \sum_{\rho=1}^lP_\mu\{F_{\rho,l}\}\bigg)\theta^{n-l} P_\mu\{H_{n-l}\le a\}\\
&=\theta\bigg(\sum_{n=1}^\infty \theta^n\sum_{\rho=1}^n P_\mu\{F_{\rho,n}\}\bigg)
\bigg(\sum_{n=1}^\infty\theta^n P_\mu\{H_n\le a\}\bigg).
\end{align*}
Plugging this into \eqref{e:est1} and noting  $\theta\in(0,1)$,  we get
\begin{align}\label{Lap}
&\sum_{n=3}^\infty \theta^nP_\mu\{H_n\le a\}\le \sum_{n=3}^\infty\theta^nP_\mu\{Q_n=0\}+12\sum_{n=1}^\infty\sum_{\rho=1}^nP_\mu\{F_{\rho, n}\}\\
&+4P\{\xi_{1,1}\le a\}\bigg(\sum_{n=1}^\infty\sum_{\rho=1}^nP_\mu\{F_{\rho, n}\}\bigg)
\sum_{n=3}^\infty \theta^nP_\mu\{H_n\le a\}.\nonumber
\end{align}

\medskip

We now install the summability of $\{P(F_{\rho,l}); 1\le\rho\le l<\infty\}$. Let $(\rho,l)$ be fixed for a while.
We start by tracking the possible site $z\in \mathbb Z$ where the early intersection ``$S_\rho=S_l=z$'' (given in the definition
of $F_{\rho,l}$) occurs. Given that $S_0=\pm 1$ and $\tilde{S}_0=\pm 1$ under $P_\mu$,  we must have that $z\in [-2,2]$. Hence,
$$
P_\mu(F_{\rho, l})=\sum_{z\in [-2,2]}P_\mu(F_{\rho, l}\cap\{S_\rho=\tilde{S}_l=z\}).
$$
Define the hitting times
$$
T_z\buildrel \Delta\over=\min\{n\ge 1; S_n=z\}\hskip.1in\hbox{and}\hskip.1in \tilde{T}_z\buildrel \Delta\over=\min\{n\ge 1; \tilde{S}_n=z\}.
$$
By the fact that ``$S_\rho=S_l$'' is an early intersection, we have
$$
\begin{aligned}
&P_\mu(F_{\rho, l}\cap\{S_\rho=\tilde{S}_l=z\})\le P_\mu\{T_z=\rho,\hskip.05in \tilde{T}_z=l\}=P_\mu\{T_z=\rho\}P_\mu\{T_z=l\},
\end{aligned}
$$
where the last step follows from the independence between $S$ and $\tilde{S}$. In summary,
\begin{align}\label{e:P-F2}
P_\mu(F_{\rho, l})\le \sum_{z\in [-2,2]}P_\mu\{T_z=\rho\}P_\mu\{T_z=l\},\hskip.2in (1\le \rho\le l<\infty).
\end{align}
Consequently,
\begin{align}\label{e:P-F3}
\sum_{l=1}^\infty\sum_{\rho=1}^l P_\mu(F_{\rho, l})\le \sum_{z\in [-2,2]}\sum_{l=1}^\infty\sum_{\rho=1}^lP_\mu\{T_z=\rho\}P_\mu\{T_z=l\}\le \sum_{z\in [-2, 2]}1=5.
\end{align}

\medskip

Choose $a>0$ sufficiently small,  so that 
$$
P\{\xi_{1,1}\le a\}<{1\over 20}.
$$
By (\ref{Lap}) and (\ref{e:P-F3}),
\begin{align}\label{Lap-1}
\sum_{n=3}^\infty \theta^nP_\mu\{H_n\le a\}\le \Big(1-20P\{\xi_{1,1}\le a\}\Big)^{-1}
\Big\{60+\sum_{n=3}^\infty\theta^nP_\mu\{Q_n=0\}\Big\}\,.
\end{align}

\begin{lemma}\label{l-u} There is a constant $C>0$ such that
\begin{align}\label{l-u-1}
P_\mu\{Q_n=0\}\le {C\over n},\hskip.2in n=1,2,\cdots,
\end{align}
and
\begin{align}\label{l-u-2}
\max_{z\in[-2,2]}P_\mu\{T_z\ge n\}\le {C\over\sqrt{n}},\hskip.2in n=1,2,\cdots.
\end{align}
\end{lemma}

\proof Under $P_\mu$, $S_0=\pm 1$ and $\tilde{S}_0=\pm 1$. When $S_0\not =\tilde{S}_0$,
$$
\{Q_n=0\}\subset
\Big\{\min_{1\le k\le n}S_k\ge -S_0,\hskip.05in \max_{1\le k\le n}\tilde{S}_k\le -\tilde{S}_0\Big\}\bigcup\Big\{\min_{1\le k\le n}\tilde{S}_k\ge -\tilde{S}_0,\hskip.05in \max_{1\le k\le n}S_k\le -S_0\Big\}.
$$
When $S_0=\tilde{S}_0$,
$$
\{Q_n=0\}\subset
\Big\{\min_{1\le k\le n}S_k\ge S_0,\hskip.05in \max_{1\le k\le n}\tilde{S}_k\le \tilde{S}_0\Big\}\bigcup\Big\{\min_{1\le k\le n}\tilde{S}_k\ge \tilde{S}_0,\hskip.05in \max_{1\le k\le n}S_k\le S_0\Big\}.
$$

By sub-additivity, symmetry and independence, we get
$$
\begin{aligned}
&P_\mu\{Q_n=0\}\\
&\le P_1\Big\{\min_{1\le k\le n}S_k\ge -1\Big\}P_{-1}\Big\{\max_{1\le k\le n}S_k\le 1\Big\}
+P_1\Big\{\min_{1\le k\le n}S_k\ge 1\Big\}P_1\Big\{\max_{1\le k\le n}S_k\le 1\Big\}\\
&=\bigg(P\Big\{\min_{1\le k\le n}S_k\ge -2\Big\}\bigg)^2+\bigg(P\Big\{\min_{1\le k\le n}S_k\ge 0\Big\}\bigg)^2
\le 2\bigg(P\Big\{\min_{1\le k\le n}S_k\ge -2\Big\}\bigg)^2.
\end{aligned}
$$

As for the hitting time, we have
$$
\begin{aligned}
&P_\mu\{T_z\ge n\}={1\over 2}P\{T_{z-1}\ge n\}+{1\over 2}P\{T_{z+1}\ge n\}\\
&={1\over 2}P\{T_{-\vert z-1\vert}\ge n\}+{1\over 2}P\{T_{-\vert z+1\vert}\ge n\}\\
&\le {1\over 2}P\Big\{\min_{1\le k\le n}S_k\ge -\vert z-1\vert\Big\}+{1\over 2}P\Big\{\min_{1\le k\le n}S_k\ge -\vert z+1\vert\Big\},
\end{aligned}
$$
where the second equality follows from symmetry.

Therefore, both (\ref{l-u-1}) and (\ref{l-u-2}) are reduced to the proof for the bound of the form
$$
P\Big\{\min_{1\le k\le n}S_k\ge -z\Big\}\le {C_z\over\sqrt{n}}\,.
$$
The bound like this looks classic and must exist somewhere in literature. For the
reader's convenience, we give a short proof here.
Since the probability on the left is non-decreasing in $z$, we may assume that $z\ge 2$. By our construction of the simple
random walk, $\vert B_t-S_k\vert\le 1$ for $\tau_k\le t\le \tau_{k+1}$. Therefore,
$$
\begin{aligned}
&P\Big\{\min_{1\le k\le n}S_k\ge -z\Big\}\le P\Big\{\min_{1\le s\le \tau_n}B_s\ge -(z-1)\Big\}\\
&\le P\Big\{\min_{1\le s\le (1-\delta)n\E\tau_1}B_s\ge -(z-1)\Big\}+P\{\tau_n\le (1-\delta) n\E\tau_1\}\\
&\le P\Big\{\vert B_{(1-\delta)n\E\tau_1}\vert\le z-1\Big\}+{1\over n\delta^2}\E\tau_1^2\le {C_z\over\sqrt{n}}
\end{aligned},
$$
where the the third step follows from Lemma \ref{lem:reflection} and the inequality (\ref{Chebyshev}). \qed

\medskip

Bringing (\ref{l-u-1}) to the bound (\ref{Lap-1}) we have 
$$
\sum_{n=3}^\infty \theta^nP_\mu\{H_n\le a\}\le C\bigg\{1+\sum_{n=1}^\infty{\theta^n\over n}\bigg\}
= C\bigg\{1+\log{1\over 1-\theta}\bigg\}.
$$
Therefore (with possibly different $C>0$),
$$
n\theta^nP_\mu\{H_n\le a\}\le\sum_{k=1}^n\theta^kP\{H_k\le k\}
\le C\log{1\over 1-\theta} \hskip0.2in (\theta\to 1^-).
$$
Taking $\theta=1-{1\over n}$ for large $n$, we have
$$
n\Big(1-{1\over n}\Big)^nP_\mu\{H_n\le a\}\le C\log n.
$$
So we have the bound
\begin{align}\label{e:ub-Hn}
  P_\mu\{H_n\le a\}=O\Big({\log n\over n}\Big)\hskip.2in (n\to\infty).
  \end{align}
This bound is weaker than the desired estimate 
\eqref{e:Hn} and therefore needs to be strengthened.
\medskip

Under the re-index $l\mapsto n-l-1$ in the second term on the
right hand side,  \eqref{e:key-estimate} becomes
$$
\begin{aligned}
  &P_\mu\{H_n\le a\}\le P_\mu\{Q_n=0\}+4P\{\xi_{1,1}\le a\}
    \sum_{l=1}^{n-2}\bigg(\sum_{\rho=1}^{n-1-l}P_\mu\{F_{\rho,n-l-1}\}\bigg)
    P_\mu\{H_l\le a\}\\
&+2P\{\xi_{1,1}\le a\}\bigg(\sum_{\rho=1}^{n-1}P_\mu\{F_{\rho,n-1}\}\bigg)+2P\{\xi_{1,1}\le a\}\bigg(\sum_{\rho=1}^{n}P_\mu\{F_{\rho,n}\}\bigg).
\end{aligned}
$$
For large integer $k$, therefore,
\begin{equation}\label{e:Hn2k}
\begin{aligned}
  &\sum_{n=2^k+1}^{2^{k+1}}P_\mu\{H_n\le a\} \le
    \sum_{n=2^k+1}^{2^{k+1}}P_\mu\{Q_n=0\}+2
    P\{\xi_{1,1}\le a\}\bigg(\sum_{n=2^k+1}^{2^{k+1}}
    \sum_{\rho=1}^{n-1}P_\mu\{F_{\rho,n-1}\}\bigg)\\
  &+2P\{\xi_{1,1}\le a\}
    \bigg(\sum_{n=2^k+1}^{2^{k+1}}\sum_{\rho=1}^{n}P_\mu\{F_{\rho,n}\}\bigg)
  \\
& +4P\{\xi_{1,1}\le a\}\sum_{n=2^k+1}^{2^{k+1}}
  \sum_{l=1}^{n-2}\bigg(\sum_{\rho=1}^{n-1-l}P_\mu\{F_{\rho,n-l-1}\}\bigg)
  P_\mu\{H_l\le a\}\\
  &\le 20+\sum_{n=2^k+1}^{2^{k+1}}P_\mu\{Q_n=0\}\\
  &+4P\{\xi_{1,1}\le a\}\sum_{n=2^k+1}^{2^{k+1}}
  \sum_{l=1}^{n-2}\bigg(\sum_{\rho=1}^{n-1-l}P_\mu\{F_{\rho,n-l-1}\}\bigg)
  P_\mu\{H_l\le a\},
\end{aligned}
\end{equation}
where the last step follows from (\ref{e:P-F3}).

For the summation in the last term on the right-hand side,
\begin{equation}\label{e:Hn2k+1}
\begin{aligned}
&\sum_{n=2^k+1}^{2^{k+1}}
\sum_{l=1}^{n-2}\bigg(\sum_{\rho=1}^{n-1-l}P_\mu\{F_{\rho,n-l-1}\}\bigg)P_\mu\{H_l\le a\}\\
&\le \sum_{l=1}^{2^{k+1}}P_\mu\{H_l\le a\}\sum_{n=\max\{l+1, 2^k\}+1}^{2^{k+1}}\sum_{\rho=1}^{n-1-l}P_\mu\{F_{\rho,n-l-1}\}\\
&\le\sum_{l=1}^{2^{k}}P_\mu\{H_l\le a\}
  \sum_{n=2^k+1}^{2^{k+1}}\sum_{\rho=1}^{n-1-l}P_\mu\{F_{\rho,n-l-1}\}\\
  &+\sum_{l=2^k+1}^{2^{k+1}}
  P_\mu\{H_l\le a\}\sum_{n=l+2}^{2^{k+1}}\sum_{\rho=1}^{n-1-l}
  P_\mu\{F_{\rho,n-l-1}\}.
\end{aligned}
\end{equation}
For the first term on the right-hand side, by (\ref{e:P-F2}) and
(\ref{e:ub-Hn}),
$$
\begin{aligned}
  &\sum_{l=1}^{2^{k}}P_\mu\{H_l\le a\}
  \sum_{n=2^k+1}^{2^{k+1}}\sum_{\rho=1}^{n-1-l}P_\mu\{F_{\rho,n-l-1}\}\\
&\le C\sum_{l=1}^{2^{k}}{\log l\over l}
  \sum_{n=2^k+1}^{2^{k+1}}\sum_{\rho=1}^{n-1-l}\sum_{z\in[-2,2]}P_\mu\{T_z=\rho\}
  P_\mu\{T_z=n-1-l\}\\
  &\le C\sum_{z\in[-2,2]}\sum_{l=1}^{2^{k}}{\log l\over l}\sum_{n=2^k+1}^\infty
    P_\mu\{T_z=n-1-l\}\\
  &\le C\sum_{z\in[-2,2]}\sum_{l=1}^{2^{k}}{\log l\over l}P_\mu\{T_z\ge 2^k-l\}\\
  &\le C\sum_{l=1}^{2^{k}}{\log l\over l}{1\over\sqrt{2^k-l+1}},
    \end{aligned}
  $$
  where the last step follows from (\ref{l-u-2}), Lemma \ref{l-u}.
  We claim that the summation on the right-hand side is bounded in $k$.
  Indeed,
  $$
  \begin{aligned}
    &\sum_{l=1}^{2^{k}}{\log l\over l}{1\over\sqrt{2^k-l+1}}
  \le 2^{-{k-1\over 2}}\sum_{l=1}^{2^{k-1}}{\log l\over l}+
  {\log 2^{k-1}\over 2^{k-1}}
      \sum_{l=2^{k-1}+1}^{2^{k}}{1\over\sqrt{2^k-l+1}}\\
    &\le C\Big\{2^{-{k-1\over 2}}(\log 2^{k-1})^2+{\log 2^{k-1}\over 2^{k-1}}2^{k-1\over 2}\Big\},
  \end{aligned}
  $$
and clearly the right-hand side is bounded in $k$.

As for the second term on the right-hand side of (\ref{e:Hn2k+1}), we use the easy bound
$$
\sum_{n=l+2}^{2^{k+1}}\sum_{\rho=1}^{n-1-l}
P_\mu\{F_{\rho,n-l-1}\}\le\sum_{n=1}^\infty\sum_{\rho=1}^nP_\mu\{F_{\rho, n}\}
\le 5,
$$
where the last step follows from (\ref{e:P-F3}).

Combining the steps after (\ref{e:Hn2k+1}), we get
$$
\sum_{n=2^k+1}^{2^{k+1}}
\sum_{l=1}^{n-2}\bigg(\sum_{\rho=1}^{n-1-l}
P_\mu\{F_{\rho,n-l-1}\}\bigg)P_\mu\{H_l\le a\}
\le C+5\sum_{n=2^k+1}^{2^{k+1}}P_\mu\{H_n\le a\}.
$$
Plugging this into \eqref{e:Hn2k}, we have
$$
\sum_{n=2^k+1}^{2^{k+1}}P_\mu\{H_n\le a\} 
\le C+\sum_{n=2^k+1}^{2^{k+1}}P_\mu\{Q_n=0\}
+20P\{\xi_{1,1}\le a\}\sum_{n=2^k+1}^{2^{k+1}}P_\mu\{H_n\le a\},
$$
or equivalently,
\begin{equation}\label{e:2k1}
\begin{aligned}
&\sum_{n=2^k+1}^{2^{k+1}}P_\mu\{H_n\le a\}
  \le\Big(1-20P\{\xi_{1,1}\le a\}\Big)^{-1}
  \Big\{C+\sum_{n=2^k+1}^{2^{k+1}}P_\mu\{Q_n=0\}\Big\}\,.
\end{aligned}
\end{equation}

In view of (\ref{l-u-1}) in Lemma \ref{l-u}, the sequence
$$
\sum_{n=2^k+1}^{2^{k+1}}P_\mu\{Q_n=0\}
\le \sum_{n=2^k+1}^{2^{k+1}} \frac Cn\le C \left(\log 2^{k+1} -\log 2^k\right)\le C
$$
is uniformly  bounded in $k$. Therefore, 
 \eqref{e:2k1} yields that  there exists a constant $C$  independent of $k$ such that 
$$
\sum_{n=2^k+1}^{2^{k+1}}P_\mu\{H_n\le a\}\le C\hskip.1in \text{ for all } k=1,2,\dots, 
$$
which leads to, noting the monotonicity of  $P_\mu\{H_n\le a\}$ in $n$,
$$
P_\mu\{H_{2^{k+1}}\le a\}\le C{1\over 2^k}
\hskip.1in \text{ for all } k=1,2,\dots
$$
By the monotonicity of the random sequence $\{H_n\}$,
this leads to  (\ref{e:Hn}). \qed

\medskip

\noindent{\textbf{Acknowledgment.}} X. Chen is partially supported by the Simons Foundation (no. 585506). J. Song is partially supported by National Natural Science Foundation of China (no. 12071256)  and  NSFC Tianyuan Key Program Project (no. 12226001).

\end{document}